\documentclass{amsart}

\usepackage{graphics}
\usepackage{epsfig}
\usepackage{amsmath}
\usepackage{amssymb,fancyhdr,txfonts,pxfonts}

\usepackage{xcolor}
\usepackage{color, colortbl}
\usepackage{graphicx}
\usepackage{marvosym}

\usepackage{paralist}
\usepackage{multicol}

\usepackage{pstricks,pst-text,pst-grad,pst-node,pst-3dplot,pstricks-add,pst-poly,pst-coil} 
\usepackage{pst-fun,pst-blur} 

\usepackage{hyperref}

\newtheorem{theorem}{Theorem}
\theoremstyle{plain}

\newtheorem{definition}{Definition}

\newtheorem{notation}{Notation}

\begin{document}

\title[Four cross-ratio maps sets of Points in Desargues Affine Plane]{Four cross-ratio maps sets of Points and their Algebraic Structures in a line on Desargues Affine Plane}

\author[Orgest ZAKA]{Orgest ZAKA}
\address{Orgest ZAKA: Department of Mathematics-Informatics, Faculty of Economy and Agribusiness, Agricultural University of Tirana, Tirana, Albania}
\email{ozaka@ubt.edu.al, gertizaka@yahoo.com}


\begin{abstract}
This paper introduces advances in the geometry of the transforms for cross ratio of four points in a line in the Desargues affine plane. The results given here have a clean, based Desargues affine plan axiomatic and definitions of addition and multiplication of points on a line in this plane, and for skew field properties. 
In this paper are studied, four types of cross-ratio maps sets of points, we discussed about for each of the 4-points of cross-ratio and we will examine the algebraic properties for each case. We are constructing four cross-ratio maps sets $\mathcal{R}^{A}_4=\left\{c_r(X,B;C,D) | \quad \forall X \in \ell^{OI} \right\}$, 
$\mathcal{R}^{B}_4=\left\{c_r(A,X;C,D) | \quad \forall X \in \ell^{OI} \right\}$,
$\mathcal{R}^{C}_4=\left\{c_r(A,B;X,D) | \quad \forall X \in \ell^{OI} \right\}$ and
$\mathcal{R}^{D}_4=\left\{c_r(A,B;C,X) | \quad \forall X \in \ell^{OI} \right\}$.
We disuse and examine algebraic properties for each case, related to the actions of addition and multiplication of points in $\ell^{OI}$ line in Desargues affine planes, which are produced by these map sets.

\end{abstract}


\maketitle

\section{Introduction and Preliminaries}
In this paper we will discuss some types of cross-ratio maps sets of points, we discussed about each of the 4-points of cross-ratio and we will examine the algebraic properties for each case.
This paper comes as a result of the influence of the recently achieved results of my research which mainly is focused on Desargues affine planes and his reports on algebraic structures, as is evident from my earlier works, see \cite{ZakaThesisPhd}, \cite{ZakaCollineations},\cite{ZakaVertex},\cite{ZakaDilauto},\cite{ZakaFilipi2016},\cite{ZakaMohammedEndo},\cite{ZakaMohammedSF}, \cite{ZakaPetersIso}, \cite{ZakaPetersOrder}. Also, in our earlier research, we think that we have made a detailed presentation about the ratio of 2 and 3 points, accompanied by the division of transformations related to these ratios, invariant and preserving transforms, we also made a presentation of the Dyck free Group and Dyck Fundamental Group, see \cite{ZakaPeters2022DyckFreeGroup}, \cite{PetersZakaFundamentalGroup2023},\cite{ZakaPeters2022InvariantPreserving}. We have also made a detailed presentation of the cross-ratio and transforms related to it, and we have presented a division of some transformations as invariant and preserving transforms, see \cite{ZakaPeters2025CrosRatio}, \cite{ZakaPeters2024}. In this paper we will discuss in turn, four cross-ratio map sets of points in Desargues affine planes.
The foundations for the study of the connections between axiomatic geometry and algebraic structures were set forth by D. Hilbert \cite{Hilbert1959geometry}. And some classic for this are,  E. Artin \cite{Artin1957GeometricAlgebra}, D.R. Huges and F.C. Piper ~\cite{HugesPiper}, H. S. M Coxeter ~\cite{CoxterIG1969}. 
Marcel Berger in \cite{Berger2009geometry12}, Robin Hartshorne
 in \cite{Hartshorne1967Foundations}, etc. 

In this paper we work in Desargues affine plane, and we study the sets of points in a line of this plane, for this reason, we will follow and remain correct to the well-known definitions and results for Desargues affine planes, presented in our previous papers, and the results for Desargues affine planes, and skew-fields presented in \cite{Pickert1973PlayfairAxiom}, \cite{Kryftis2015thesis}, \cite{Artin1957GeometricAlgebra}, \cite{HugesPiper}, \cite{Hilbert1959geometry}, \cite{Cohn2008sf}, \cite{Herstein1968NR},\cite{Rotman2015AMAlgebra},\cite{Lam2001GTMalgebra}.\\
\textbf{Desargues Affine Plane:} 
Let $\mathcal{P}$ be a nonempty space, $\mathcal{L}$ a nonempty subset of $\mathcal{P}$. The elements $P$ of $\mathcal{P}$ are points and an element $\ell$ of $\mathcal{L}$ is a line. 
\begin{definition}
The incidence structure $\mathcal{A}=(\mathcal{P}, \mathcal{L},\mathcal{I})$, called affine plane, where satisfies the above axioms:
\begin{description}
	\item[1$^o$] For each points $\left\{P,Q\right\}\in \mathcal{P}$, there is exactly one line $\ell\in \mathcal{L}$ such that $\left\{P,Q\right\}\in \ell$.
\item[2$^o$] For each point $P\in \mathcal{P}, \ell\in \mathcal{L}, P \not\in \ell$, there is exactly one line $\ell'\in \mathcal{L}$ such that
$P\in \ell'$ and $\ell\cap \ell' = \emptyset$\ (Playfair Parallel Axiom~\cite{Pickert1973PlayfairAxiom}).   Put another way,
if the point $P\not\in \ell$, then there is a unique line $\ell'$ on $P$ missing $\ell$~\cite{Prazmowska2004DemoMathDesparguesAxiom}.
\item[3$^o$] There is a 3-subset of points $\left\{P,Q,R\right\}\in \mathcal{P}$, which is not a subset of any $\ell$ in the plane.   Put another way, there exist three non-collinear points $\mathcal{P}$~\cite{Prazmowska2004DemoMathDesparguesAxiom}.
\end{description}
\end{definition}

\begin{figure}[htbp]
	\centering
		\includegraphics[width=0.85\textwidth]{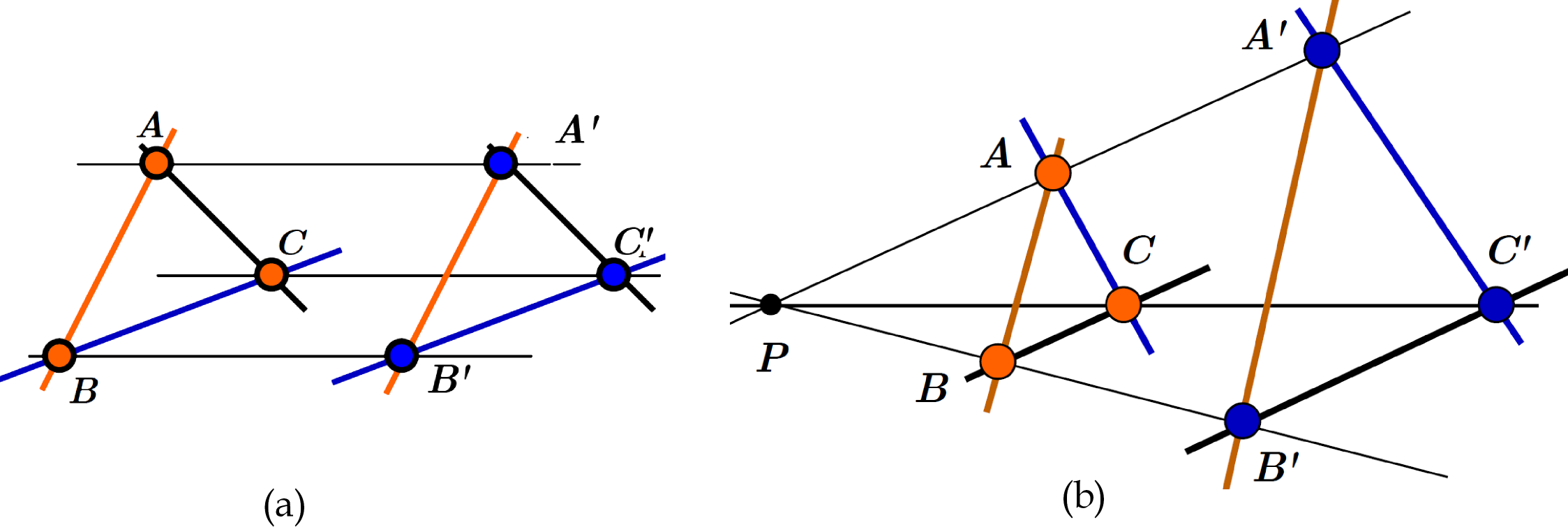}
	\caption{Desargues Axioms: (a) For parallel lines $\ell^{AA'} \parallel \ell^{BB'} \parallel \ell^{CC'}$; (b) For lines which are cutting in a single point $P$,  $\ell^{AA'} \cap \ell^{BB'} \cap \ell^{CC'}=P$.}
		\label{fig:DesarguesAxiom}
\end{figure}

\emph{\bf Desargues' Axiom, circa 1630}~\cite[\S 3.9, pp. 60-61] {Kryftis2015thesis}~\cite{Szmielew1981DesarguesAxiom}.   Let $A,B,C,A',B',C'\in \mathcal{P}$ and let pairwise distinct lines  $\ell^{AA'} , \ell^{BB'}, \ell^{CC'}, \ell^{AC}, \ell^{A'C'}\in \mathcal{L}$ such that
\begin{align*}
\ell^{AA'} \parallel \ell^{BB'} \parallel \ell^{CC'} \ \mbox{(Fig.~\ref{fig:DesarguesAxiom}(a))} &\ \mbox{\textbf{or}}\
\ell^{AA'} \cap \ell^{BB'} \cap \ell^{CC'}=P.
 \mbox{(Fig.~\ref{fig:DesarguesAxiom}(b) )}\\
 \mbox{and}\  \ell^{AB}\parallel \ell^{A'B'}\ &\ \mbox{and}\ \ell^{BC}\parallel \ell^{B'C'}.\\
A,B\in \ell^{AB}, A'B'\in \ell^{A'B'}, A,C\in \ell^{AC}, &\  \mbox{and}\ A'C'\in \ell^{A'C'}, B,C\in \ell^{BC},  B'C'\in \ell^{B'C'}.\\
A\neq C, A'\neq C', &\ \mbox{and}\ \ell^{AB}\neq \ell^{A'B'}, \ell^{BC}\neq \ell^{B'C'}.
\end{align*}
Then $\boldsymbol{\ell^{AC}\parallel \ell^{A'C'}}$.  

\noindent A {\bf Desargues affine plane} is an affine plane that satisfies Desargues' Axiom. 
\begin{notation}
Three vertexes $ABC$ and $A'B'C'$, which, fulfilling the conditions of the Desargues Axiom, we call \emph{'Desarguesian'}.
\end{notation}

\begin{figure}[htbp]
\centering%
\includegraphics[width=0.92\textwidth]{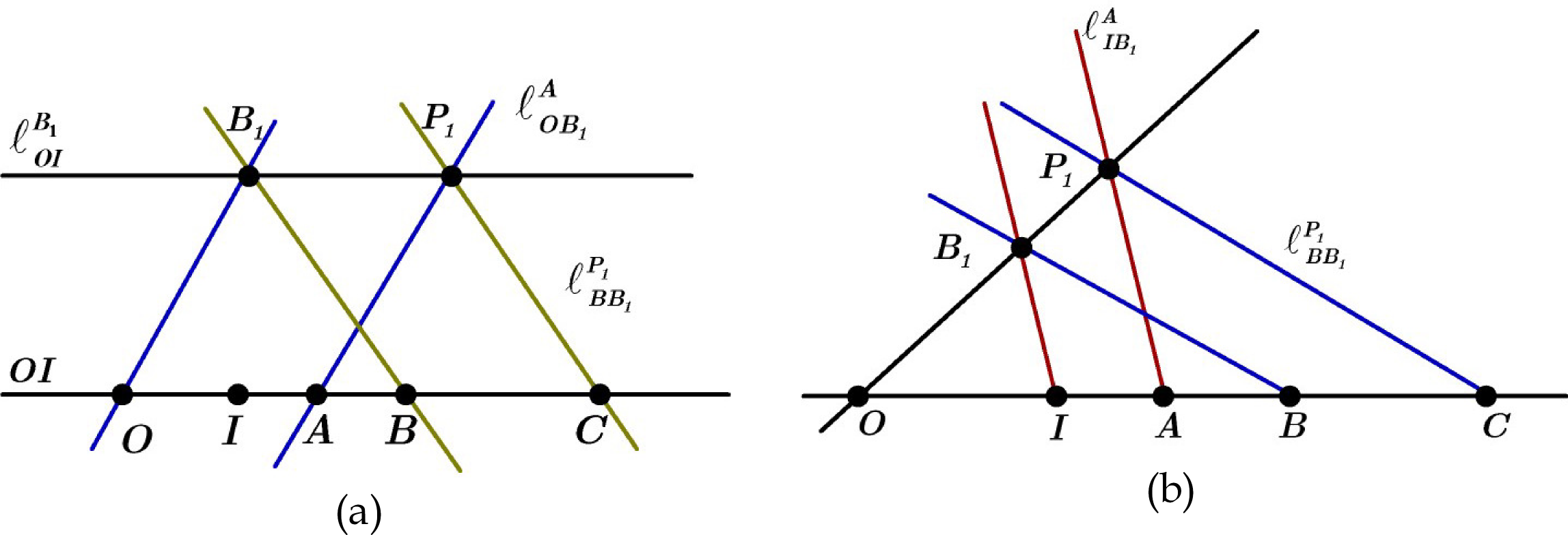}
\caption{ (a) Addition of points in a line in affine plane, 
(b) Multiplication of points in a line in affine plane}
\label{fig:FigureAdMult}
\end{figure}

\textbf{Addition and Multiplication of points in a line of Desargues affine plane \cite{ZakaThesisPhd}:}
The process of construct the points $C$ for addition (Figure \ref{fig:FigureAdMult} (a)) and multiplication (Figure \ref{fig:FigureAdMult} (b)) of points in $\ell^{OI}-$line in affine plane, is presented in the tow algorithmic form  
\begin{multicols}{2}
\textsc{Addition Algorithm}
\begin{description}
	\item[Step.1] $B_{1}\notin \ell^{OI}$
	\item[Step.2] $\ell_{OI}^{B_{1}}\cap \ell_{OB_{1}}^{A}=P_{1}$
	\item[Step.3] $\ell_{BB_{1}}^{P_{1}}\cap \ell^{OI}=C(=A+B)$
\end{description}

\textsc{Multiplication Algorithm}
\begin{description}
	\item[Step.1] $B_{1}\notin \ell^{OI}$
	\item[Step.2] $\ell_{IB_{1}}^{A}\cap \ell^{OB_{1}}=P_{1}$
	\item[Step.3] $\ell_{BB_{1}}^{P_{1}}\cap \ell^{OI}=C(=A\cdot B)$
\end{description}
\end{multicols}
In \cite{ZakaThesisPhd} and \cite{FilipiZakaJusufi}, we have prove that $(\ell^{OI}, +, \cdot)$ is a skew field in Desargues affine plane, and is field (commutative skew field) in the Papus affine plane.
\begin{definition} \label{ratio2points}
\cite{ZakaPeters2022DyckFreeGroup} Lets have two different points $A,B \in \ell^{OI}-$line, and $B\neq O$, in Desargues affine plane. We define as ratio of this tow points, a point $R\in \ell^{OI}$, such that,
\[R=B^{-1}A, \qquad \text{ we mark this, with,} \qquad 
R=r(A:B)=B^{-1}A
\]
\end{definition}
For a 'ratio-point' $R \in \ell^{OI}$, and for point $B\neq O$ in line $\ell^{OI}$, is a unique defined point, $A \in \ell^{OI}$, such that $R=B^{-1}A=r(A:B)$.
\begin{figure}[htbp]
	\centering
		\includegraphics[width=0.75\textwidth]{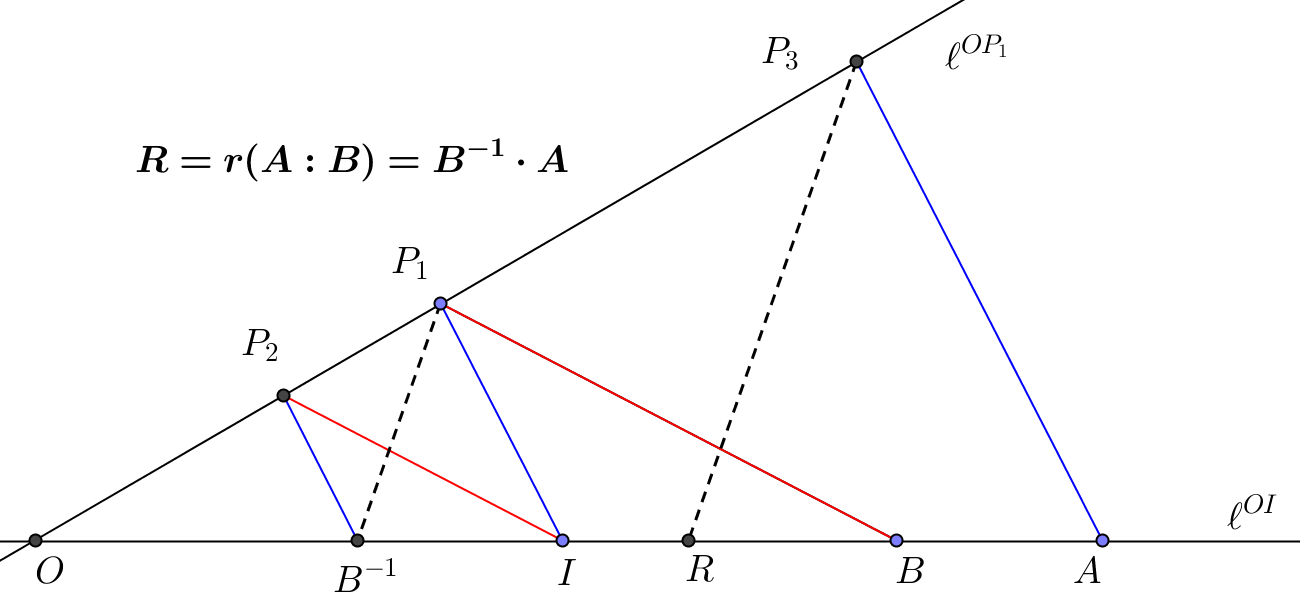}
	\caption{Ilustrate the Ratio-Point, of 2-Points in a line of Desargues affine plane $R=r(A:B)=B^{-1}A$.}
	\label{Ratio2points}
\end{figure}
In paper \cite{ZakaPeters2022DyckFreeGroup} we have presented a detailed study and more results related to the ratio of 3 points.
\begin{definition}\label{ratiodef}
If $A, B, C$ are three points on a line $\ell^{OI}$ (collinear) in Desargues affine plane, then we define their \textbf{ratio} to be a point $R \in \ell^{OI}$, such that:
\[
(B-C)\cdot R=A-C, \quad \mbox{concisely}\quad R=(B-C)^{-1}(A-C),
\]
and we mark this with  $r(A,B;C)= (B-C)^{-1}(A-C)$.
\end{definition}
\begin{figure}[htbp]
	\centering
		\includegraphics[width=0.9\textwidth]{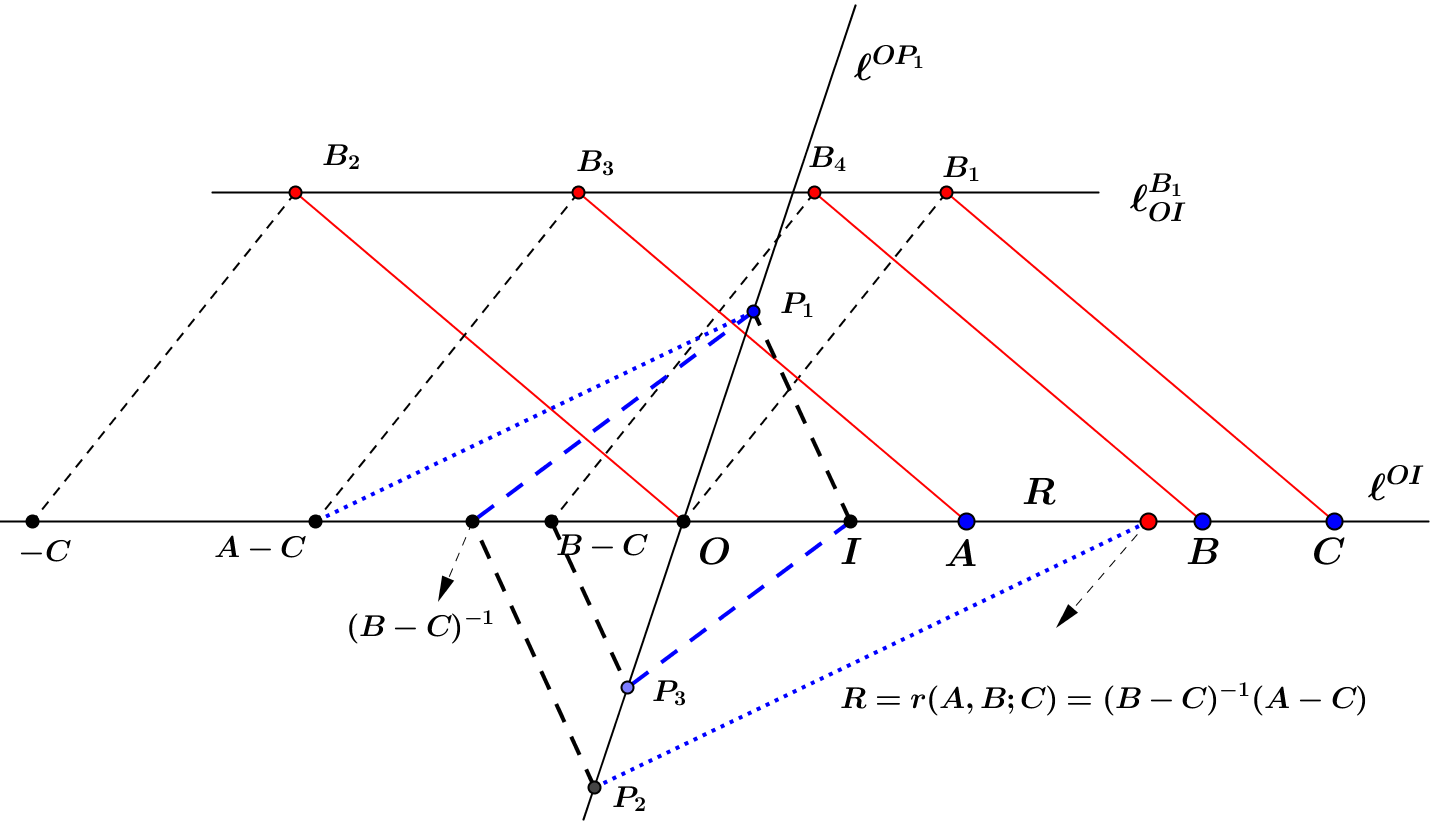}
	\caption{Ratio of 3-Points in a line of Desargues affine plane $R=r(A,B;C)$.}
		\label{ratio3points1}
\end{figure}
\textbf{Cross-Ratio in a line of Desargues affine plane:} 
Let us have the line $\ell^{OI}$ in Desarges affine plane $\mathcal{A_{D}}$, and four points, $A, B, C, D \in \ell^{OI}$
\begin{definition}\label{cross-ratio.def}
If $A, B, C, D$ are four points on a line $\ell^{OI}$ in Desarges affine plane $\mathcal{A_{D}}$, no three of them equal, then we define their cross ratio to be a point:
\[c_r(A,B;C,D)=\left[(A-D)^{-1}(B-D)\right]\left[(B-C)^{-1}(A-C)\right]
\]
\end{definition}

\section{Four 'cross-ratio map sets' of points in a line of Desargues affine plane and their algebraic properties.}

In this section we will discuss four types of cross-ratio maps sets of points, we discussed about for each of the 4-points of cross-ratio and we will examine the algebraic properties for each case. We are constructing four cross-ratio maps sets, below.
\begin{enumerate}
	\item  $\mathcal{R}^{A}_4=\left\{c_r(X,B;C,D) | \quad \forall X \in \ell^{OI} \right\}$,  where $B,C,D\in \ell^{OI}$ ;
	\item $\mathcal{R}^{B}_4=\left\{c_r(A,X;C,D) | \quad \forall X \in \ell^{OI} \right\}$ where $A,C,D\in \ell^{OI}$;
	\item $\mathcal{R}^{C}_4=\left\{c_r(A,B;X,D) | \quad \forall X \in \ell^{OI} \right\}$ where$A,B,D \in \ell^{OI}$;
	\item $\mathcal{R}^{D}_4=\left\{c_r(A,B;C,X) | \quad \forall X \in \ell^{OI} \right\}$ where $A,B,C \in \ell^{OI}$.
\end{enumerate}
We also note that the points of a line on Desargues affine plane, are closed under two operations: "Addition" and "Multiplication" of points in this line.
\subsection{Algebraic properties for the 'cross-ratio map sets' $\mathcal{R}^{A}_4$, set of Mobi\"us points} $  $\\
We see that the first, $\mathcal{R}^{A}_4$, is the 'cross-ratio map sets' of Mobi\"us points, we are calling this "set of Mobi\"us points" after the Mobi\"us transformation for cross-ratio points i n Desargues affine plane, which we studied in the paper \cite{ZakaPeters2024}, we are introducing this definition,
\begin{definition} \label{mobius.transform.def}
Lets have three "bass points", $B,C,D\in \ell^{OI}.$ Mobi\"us transform for cross-ratio, called the map, $ \mu: \ell^{OI} \to \ell^{OI}, $ which satisfies the condition, 
\[
\forall X \in \ell^{OI}, \quad \mu(X)=c_r(X,B;C,D).
\]
\end{definition}
From Mobi\"us Transform in definition \ref{mobius.transform.def}, but also by the ratio map set for the ratio of two and three points in a line in Desargues affine plane, we have constructed the set $\mathcal{R}^{A}_4$, which is the set,
\[ 
\mathcal{R}^{A}_4=\left\{c_r(X,B;C,D) | \quad \forall X \in \ell^{OI} \right\}
\]
We call this set also the set of \emph{'Mobi\"us points'} in $\ell^{OI}-$line on Desargues affine plane.

From Mobi\"us transform definition \ref{mobius.transform.def} we have that
\[\mu(X)=c_r(X,B;C,D)=[(X-D)^{-1}(B-D)][(B-C)^{-1}(X-C)]
\]
so, for elements (Mobi\"us points) of $\mathcal{R}_4$ we market with $\mu(X)$.

We will see in the following Theorems, the properties that are fulfilled related to the addition and multiplication of 'Mobi\"us points' in $\ell^{OI}-$line, so, we will see, the algebraic structures that are created, related to the addition and multiplication of 'Mobi\"us points' in $\ell^{OI}-$line in Desargues affine plane.

\begin{theorem} \label{thm.11}
Let's have three different-fixed points $B, C, D \in \ell^{OI}$  and this points are different for point $O$, in  $\ell^{OI}-$line. The Mobi\"us transforms $\mu(\cdot)$ from $\mathcal{R}_4$, with addition of points, fulfills the properties: \\
(1) is \emph{Associative},\\
(2) is commutative and, \\
(3) it has a neutral element (zero element).
\end{theorem}
\proof
\textbf{1.} (Associativity) $\forall \mu(X),\mu(Y),\mu(Z) \in \mathcal{R}^{A}_4$, we have three points $ X,Y,Z \in \ell^{OI}$, and, we have that:
\[
\begin{aligned}
\left[\mu(X)+\mu(Y)\right]+\mu(Z)
&= \{ [(X-D)^{-1}(B-D)][(B-C)^{-1}(X-C)] \\
&+[(Y-D)^{-1}(B-D)][(B-C)^{-1}(Y-C)] \} \\
&+ [(Z-D)^{-1}(B-D)][(B-C)^{-1}(Z-C)] 
\end{aligned}
\]
cross-ratio points are point on $\ell^{OI}-$line, and 
from asociatiove properties of addition groups $(\ell^{OI},+)$, so we have that,
\[
\begin{aligned}
\left[\mu(X)+\mu(Y)\right]+\mu(Z)
&= [(X-D)^{-1}(B-D)][(B-C)^{-1}(X-C)] \\
&+\{ [(Y-D)^{-1}(B-D)][(B-C)^{-1}(Y-C)]  \\
&+ [(Z-D)^{-1}(B-D)][(B-C)^{-1}(Z-C)] \} \\
&=\mu(X)+\left[\mu(Y)+\mu(Z)\right].
\end{aligned}
\]
Hence,
\[
[\mu(X)+\mu(Y)]+\mu(Z)=\mu(X)+[\mu(Y)+\mu(Z)], \forall X,Y,Z \in \ell^{OI}.
\]
\textbf{2.} (Commuativity) $\forall X,Y \in \ell^{OI}$, we have that, the 'cross-ratio points' are point on $\ell^{OI}-$line, and 
from commutative properties of addition groups $(\ell^{OI},+)$, we have that,
\[
\begin{aligned}
\mu(X)+\mu(Y)&=c_r(X,B;C,D)+c_r(Y,B;C,D)\\
&=[(X-D)^{-1}(B-D)][(B-C)^{-1}(X-C)] +[(Y-D)^{-1}(B-D)][(B-C)^{-1}(Y-C)]\\
&=[(Y-D)^{-1}(B-D)][(B-C)^{-1}(Y-C)]+[(X-D)^{-1}(B-D)][(B-C)^{-1}(X-C)] \\
&=c_r(Y,B;C,D)+c_r(X,B;C,D)\\
&=\mu(Y)+\mu(X),
\end{aligned}
\]
so
\[ \mu(X)+\mu(Y)=\mu(Y)+\mu(X),  \forall X,Y \in \ell^{OI}\]

\textbf{3.} ('zero-element') $\exists Z \in \ell^{OI}, \forall X \in \ell^{OI}$, we have that:
\[ c_r(X,B;C,D)+c_r(Z,B;C,D)=c_r(X,B;C,D) \]
cross-ratio point are point of $\ell^{OI}$-line, and '$+$' is addition of points in $\ell^{OI}$-line, we have prove that the 'zero-element' for addition of points in $\ell^{OI}$-line, is point $O$, thus, have that,
\[ c_r(Z,B;C,D)=O \]
For point $O$ (zero element for addition of points), have
\[
\begin{aligned}
c_r(Z,B;C,D)&=O\\
[(Z-D)^{-1}(B-D)][(B-C)^{-1}(Z-C)] &=O\\
\end{aligned}
\]
also we have that the skew field $(\ell^{OI},+, \cdot )$ do not have 'divisors of zero', therefore,
 \[ [(Z-D)^{-1}(B-D)]=O \quad\text{or}\quad [(B-C)^{-1}(Z-C)]=O \]
more,

 \[ (Z-D)^{-1}=O \quad\text{or}\quad (B-D)=O, \quad\text{or}\quad  (B-C)^{-1}=O, \quad\text{or}\quad (Z-C)=O .\]

But, we suppose that the point $B \neq D$,  $B \neq C$ and $D \neq C$, therefore we have that,
\[
B-D \neq O \quad\text{and}\quad O\neq (B-C)^{-1} \neq \infty
\]
so have that, $(Z-D)^{-1}=O \quad\text{or}\quad (Z-C)=O$.

Hence, 
\[Z-C=O \Leftrightarrow Z=C
\]
so 'zero-element' is $\mu(C)$.
\qed

\begin{theorem}\label{thm.12}
Let's have three different points $B,C,D$ and different for point $O$, in  $\ell^{OI}-$line. The cross-ratio map set (set of 'Mobi\"us points'), with multiplication of points in a line of Desargues affine plane $\left( \mathcal{R}_4, \cdot \right)$ is  Group.
\end{theorem}
\proof
\textbf{1.} (Associativity) $\forall X,Y,Z \in \ell^{OI}$, we have that:
\[
\begin{aligned}
\left[\mu(X) \cdot \mu(Y)\right]\cdot \mu(Z)&=
\{[(X-D)^{-1}(B-D)][(B-C)^{-1}(X-C)] \\
& \cdot [(Y-D)^{-1}(B-D)][(B-C)^{-1}(Y-C)] \} \\
& \cdot [(Z-D)^{-1}(B-D)][(B-C)^{-1}(Z-C)]
\end{aligned}
\]
all the factors are points of the line $\ell^{OI}$, 
they are also elements of skew fields, therefore from the
association property, we can move the brackets, so we have that,
\[
\begin{aligned}
\left[\mu(X) \cdot \mu(Y)\right]\cdot \mu(Z)&=
[(X-D)^{-1}(B-D)][(B-C)^{-1}(X-C)] \\
& \cdot \{ [(Y-D)^{-1}(B-D)][(B-C)^{-1}(Y-C)]  \\
& \cdot [(Z-D)^{-1}(B-D)][(B-C)^{-1}(Z-C)] \}\\
&=\mu(X) \cdot [\mu(Y) \cdot \mu(Z) ]
\end{aligned}
\]
Hence,
\[
\left[\mu(X) \cdot \mu(Y)\right]\cdot \mu(Z)=\mu(X) \cdot [\mu(Y) \cdot \mu(Z) ]
\]
\textbf{2.} ('unitary') $\exists E \in \ell^{OI}, \forall X \in \ell^{OI}$, we have that:
\[
\mu(X)\cdot \mu(E)=\mu(X)\]
so, 
\[
\{[(X-D)^{-1}(B-D)][(B-C)^{-1}(X-C)] \}
 \cdot \{[(E-D)^{-1}(B-D)][(B-C)^{-1}(E-C)] \}  \]
\[ =[(X-D)^{-1}(B-D)][(B-C)^{-1}(X-C)]\]
Both elements $\mu(X), \mu(E)$ are points in the line $\ell^{OI}$, we also know that the unitary element related to the multiplication of points in $\ell^{OI}-$ line is \emph{unique point} $I$, for this reason we have $\mu(E)=I$. We see that, for $E=B$, we have,
\[
\mu(E)=\mu(B)=[(B-D)^{-1}(B-D)][(B-C)^{-1}(B-C)]=I\cdot I=I.
\]
so unitary-element element is $\mu(B)$.

\textbf{3.} ('inverse') $\forall X \in \ell^{OI}, \exists Y \in \ell^{OI}$, we have that:
\[
\mu(X)\cdot \mu(Y)=\mu(B)\]
so,
\[
\mu(X)\cdot \mu(Y)=I\Rightarrow \mu(Y)=\{\mu(X) \}^{-1} \]
\[
\{[(X-D)^{-1}(B-D)][(B-C)^{-1}(X-C)] \}
 \cdot \{[(Y-D)^{-1}(B-D)][(B-C)^{-1}(Y-C)]\} \]
\[  =[(B-D)^{-1}(B-D)][(B-C)^{-1}(B-C)]=I\]
hence,
\[
[(Y-D)^{-1}(B-D)][(B-C)^{-1}(Y-C)]=\{ [(X-D)^{-1}(B-D)][(B-C)^{-1}(X-C)] \}^{-1} \]
\[
=\{[(B-C)^{-1}(X-C)] \}^{-1}\{ [(X-D)^{-1}(B-D)]\}^{-1} \]
\[\mu(Y) =[(X-C)^{-1}(B-C)][(B-D)^{-1}(X-D)] \]
Hence,
\[ \{\mu(X) \}^{-1}=c_r(X, B; D,C) \]
\qed

\begin{theorem} \label{thm.13}
In the cross-ratio-maps-set $\mathcal{R}_4=\left\{\mu(X)=c_r(X,B;C,D) | \quad \forall X \in \ell^{OI} \right\}$, for three elements $\mu(X)$, $\mu(Y)$, $\mu(Z)$ of its, have true, that equations
\begin{enumerate}
	\item $\mu(X)\cdot [\mu(Y)+\mu(Z)]=\mu(X)\cdot \mu(Y) +\mu(X)\cdot \mu(Z)$, and,
	\item $[\mu(X)+\mu(Y)]\mu(Z)=\mu(X)\cdot \mu(Z) +\mu(Y)\cdot \mu(Z)$
\end{enumerate}
\end{theorem}
\proof The points $\mu(X) , \mu(Y), \mu(Z)$ are point of $\ell^{OI}-$line, and are true this equations from fact that $(\ell^{OI}, +, \cdot)$ is a skew field.
\qed

%


\subsection{Algebraic properties for the 'cross-ratio map sets' $\mathcal{R}^{B}_4$-set of points} $  $\\

We have that,
\[ 
\mathcal{R}^{B}_4=\left\{c_r(A,X;C,D) | \quad \forall X \in \ell^{OI} \right\}
\]

We mark elements of this set with $f_B(X)$ and, we have that
\[f_B(X)=c_r(A,X;C,D)=[(A-D)^{-1}(X-D)][(X-C)^{-1}(A-C)]
\]

We will see in the following Theorems, the properties that are fulfilled related to the addition and multiplication of points' in $\ell^{OI}-$line, so, we will see, the algebraic structures that are created, related to the addition and multiplication of this points in $\ell^{OI}-$line in Desargues affine plane.

We will not present proofs for some of the algebraic properties in this part, as they are the same as the proofs in the first part. We will present evidence only in cases where we have changes.

\begin{theorem} \label{thm.21}
Let's have three different-fixed points $A, C, D \in \ell^{OI}$  and this points are different for point $O$, in  $\ell^{OI}-$line. The transforms $f_B(\cdot)$ from $\mathcal{R}^{B}_4$, with addition of points, fulfills the properties: \\
(1) Associative,\\
(2) commutative and, \\
(3) it has a neutral element (zero element).
\end{theorem}
\proof  
\textbf{1.} Associative  and \textbf{2.} Commutative properties 
prove similar to Theorem \ref{thm.11}.

\textbf{3.} ('zero-element') $\exists Z \in \ell^{OI}, \forall X \in \ell^{OI}$, we have that:
\[ c_r(A,X;C,D)+c_r(A,Z;C,D)=c_r(A,X;C,D) \]
cross-ratio point are point of $\ell^{OI}$-line, and '$+$' is addition of points in $\ell^{OI}$-line, we have prove that the 'zero-element' for addition of points in $\ell^{OI}$-line, is point $O$, thus, have that,
\[ c_r(A,Z;C,D)=O \]
For point $O$ (zero element for addition of points), have
\[
\begin{aligned}
c_r(A,Z;C,D)&=O\\
[(A-D)^{-1}(Z-D)][(Z-C)^{-1}(A-C)] &=O\\
\end{aligned}
\]
also we have that the skew field $(\ell^{OI},+, \cdot )$ do not have 'divisors of zero', therefore,
 \[ [(A-D)^{-1}(Z-D)]=O \quad\text{or}\quad [(Z-C)^{-1}(A-C)]=O \]
more,

 \[ (A-D)^{-1}=O \quad\text{or}\quad (Z-D)=O, \quad\text{or}\quad(Z-C)^{-1}=O, \quad\text{or}\quad (A-C)=O .\]

But, we suppose that the point $A \neq D$,  $A \neq C$ and $D \neq C$, therefore we have that,
\[
Z-D \neq O \quad\text{and}\quad O\neq (A-D)^{-1} \neq \infty
\]
so have that, $(Z-C)^{-1}=O \quad\text{or}\quad (Z-D)=O $. 

Hence, 
\[Z-C=O \Leftrightarrow Z=C
\]
so 'zero-element' is $f_B(C)$.
\qed

\begin{theorem}\label{thm.22}
Let's have three different points $A,C,D$ and different for point $O$, in  $\ell^{OI}-$line. The cross-ratio map set with multiplication of points in a line of Desargues affine plane $\left( \mathcal{R}^{B}_4, \cdot \right)$ is  Group.
\end{theorem}
\proof
\textbf{1.} Associative properties prove similar to Theorem \ref{thm.12} 

\textbf{2.} ('unitary') $\exists E \in \ell^{OI}, \forall X \in \ell^{OI}$, we have that:
\[
f_B(X)\cdot f_B(E)=f_B(X)\]
so, 
\[
\{[(A-D)^{-1}(X-D)][(X-C)^{-1}(A-C)] \}
 \cdot \{[(A-D)^{-1}(E-D)][(E-C)^{-1}(A-C)] \}  \]
\[ =[(A-D)^{-1}(X-D)][(X-C)^{-1}(A-C)]\]
Both elements $f_B(X), f_B(E)$ are points in the line $\ell^{OI}$, we also know that the unitary element related to the multiplication of points in $\ell^{OI}-$ line is \emph{unique point} $I$, for this reason we have $f_B(E)=I$. We see that, for $E=A$, we have,
\[
f_B(E)=f_B(A)=[(A-D)^{-1}(A-D)][(A-C)^{-1}(A-C)]=I\cdot I=I.
\]
so unitary-element element is $f_B(A)$.

\textbf{3.} ('inverse') $\forall X \in \ell^{OI}, \exists Y \in \ell^{OI}$, we have that:
\[
f_B(X)\cdot f_B(Y)=f_B(A)\]
so,
\[
f_B(X)\cdot f_B(Y)=I\Rightarrow f_B(Y)=\{f_B(X) \}^{-1} \]
\[
\{[(A-D)^{-1}(X-D)][(X-C)^{-1}(A-C)] \}
 \cdot \{[(A-D)^{-1}(Y-D)][(Y-C)^{-1}(A-C)]\} \]
\[  =[(A-D)^{-1}(A-D)][(A-C)^{-1}(A-C)]=I\]
hence,
\[
[(A-D)^{-1}(Y-D)][(Y-C)^{-1}(A-C)]=\{ [(A-D)^{-1}(X-D)][(X-C)^{-1}(A-C)] \}^{-1} \]
\[
=\{[(X-C)^{-1}(A-C)] \}^{-1}\{[(A-D)^{-1}(X-D)]\}^{-1} \]
\[f_B(Y) =[(A-C)^{-1}(X-C)][(X-D)^{-1}(A-D)] \]
Hence,
\[ \{f_B(X) \}^{-1}=c_r(A, X; D,C) \]
\qed

\begin{theorem} \label{thm.23}
In the cross-ratio-maps-set $\mathcal{R}^{B}_4=\left\{f_B(X)=c_r(A,X;C,D) | \quad \forall X \in \ell^{OI} \right\}$, for three elements $f_B(X)$, $f_B(Y)$, $f_B(Z)$ of its, have true, that equations
\begin{enumerate}
	\item $f_B(X)\cdot [f_B(Y)+f_B(Z)]=f_B(X)\cdot f_B(Y) +f_B(X)\cdot f_B(Z)$, and,
	\item $[f_B(X)+f_B(Y)]f_B(Z)=f_B(X)\cdot f_B(Z) +f_B(Y)\cdot f_B(Z)$
\end{enumerate}
\end{theorem}
\proof In the same way as Theorem \ref{thm.13}.
\qed

\subsection{Algebraic properties for the 'cross-ratio map sets' $\mathcal{R}^{C}_4$-set of points} $  $\\
We have that,
\[ 
\mathcal{R}^{C}_4=\left\{c_r(A,B;X,D) | \quad \forall X \in \ell^{OI} \right\}
\]
We mark elements of this set with $f_C(X)$ and, we have that
\[f_C(X)=c_r(A,B;X,D)=[(A-D)^{-1}(B-D)][(B-X)^{-1}(A-X)]
\]

\begin{theorem} \label{thm.31}
Let's have three different-fixed points $A, B, D \in \ell^{OI}$  and this points are different for point $O$, in  $\ell^{OI}-$line. The transforms $f_C(\cdot)$ from $\mathcal{R}^{C}_4$, with addition of points, fulfills the properties: \\
(1) Associative,\\
(2) commutative and, \\
(3) it has a neutral element (zero element).
\end{theorem}
\proof  
\textbf{1.} Associative  and \textbf{2.} Commutative properties 
prove similar to Theorem \ref{thm.11}.

\textbf{3.} ('zero-element') $\exists Z \in \ell^{OI}, \forall X \in \ell^{OI}$, we have that:
\[ c_r(A,B;X,D)+c_r(A,B;X,D)=c_r(A,B;X,D) \]
cross-ratio point are point of $\ell^{OI}$-line, and '$+$' is addition of points in $\ell^{OI}$-line, we have prove that the 'zero-element' for addition of points in $\ell^{OI}$-line, is point $O$, thus, have that,
\[ c_r(A,B;X,D)=O \]
For point $O$ (zero element for addition of points), have
\[ c_r(A,B;Z,D)=O \]
so,
\[ [(A-D)^{-1}(B-D)][(B-Z)^{-1}(A-Z)] =O \]
also we have that the skew field $(\ell^{OI},+, \cdot )$ do not have 'divisors of zero', therefore,
 \[ [(A-D)^{-1}(B-D)]=O \quad\text{or}\quad [(B-Z)^{-1}(A-Z)]=O \]
more,

 \[ (A-D)^{-1}=O \quad\text{or}\quad (B-D)=O, \quad\text{or}\quad 
(B-Z)^{-1}=O ,\quad\text{or}\quad (A-Z)=O .\]

But, we suppose that the point $A \neq B$,  $A \neq D$ and $B \neq D$, therefore we have that,
\[
B-D \neq O \quad\text{and}\quad O\neq (A-D)^{-1} \neq \infty
\]
so have that, $(B-Z)^{-1}=O \quad\text{or}\quad (A-Z)=O$ .

Hence, 
\[A-Z=O \Leftrightarrow Z=A
\]
so 'zero-element' is $f_C(A)$.
\qed

\begin{theorem}\label{thm.32}
Let's have three different points $A,C,D$ and different for point $O$, in  $\ell^{OI}-$line. The cross-ratio map set with multiplication of points in a line of Desargues affine plane $\left( \mathcal{R}^{B}_4, \cdot \right)$ is  Group.
\end{theorem}
\proof
\textbf{1.} Associative properties prove similar to Theorem \ref{thm.12}

\textbf{2.} ('unitary') $\exists E \in \ell^{OI}, \forall X \in \ell^{OI}$, we have that:
\[
f_C(X)\cdot f_C(E)=f_C(X)\]
so, 
\[
\{[(A-D)^{-1}(B-D)][(B-X)^{-1}(A-X)] \}
 \cdot \{[(A-D)^{-1}(B-D)][(B-E)^{-1}(A-E)] \}  \]
\[ =[(A-D)^{-1}(B-D)][(B-X)^{-1}(A-X)]\]
Both elements $f_C(X), f_C(E)$ are points in the line $\ell^{OI}$, we also know that the unitary element related to the multiplication of points in $\ell^{OI}-$ line is \emph{unique point} $I$, for this reason we have $f_C(E)=I$. We see that, for $E=D$, we have,
\[
\begin{aligned}
f_C(E)&=f_C(D)\\
&=[(A-D)^{-1}(B-D)][(B-D)^{-1}(A-D)]\\
&=(A-D)^{-1}[(B-D)(B-D)^{-1}](A-D)\\
&=(A-D)^{-1}[I](A-D)\\
&=(A-D)^{-1}(A-D) \\
&=I
\end{aligned}
\]
so unitary-element element is $f_C(D)$.

\textbf{3.} ('inverse') $\forall X \in \ell^{OI}, \exists Y \in \ell^{OI}$, we have that:
\[
f_C(X)\cdot f_C(Y)=I=f_C(D)\]
so,
\[
f_C(X)\cdot f_C(Y)=I\Rightarrow f_C(Y)=\{f_C(X) \}^{-1} \]
\[
\{[(A-D)^{-1}(B-D)][(B-X)^{-1}(A-X)] \}
 \cdot \{[(A-D)^{-1}(B-D)][(B-Y)^{-1}(A-Y)]\} \]
\[  =[(A-D)^{-1}(B-D)][(B-D)^{-1}(A-D)]=I\]
hence,
\[
[(A-D)^{-1}(B-D)][(B-Y)^{-1}(A-Y)]=\{ [(A-D)^{-1}(B-D)][(B-X)^{-1}(A-X)] \}^{-1} \]
\[
=\{[(B-X)^{-1}(A-X)] \}^{-1}\{[(A-D)^{-1}(B-D)]\}^{-1} \]
\[f_B(Y) =[(A-X)^{-1}(B-X)][(B-D)^{-1}(A-D)] \]
Hence,
\[ \{f_B(X) \}^{-1}=c_r(A, B; D,X) \]
\qed

\begin{theorem} \label{thm.33}
In the cross-ratio-maps-set $\mathcal{R}^{C}_4=\left\{f_C(X)=c_r(A,B;X,D) | \quad \forall X \in \ell^{OI} \right\}$, for three elements $f_C(X)$, $f_C(Y)$, $f_C(Z)$ of its, have true, that equations
\begin{enumerate}
	\item $f_C(X)\cdot [f_C(Y)+f_C(Z)]=f_C(X)\cdot f_C(Y) +f_C(X)\cdot f_B(Z)$, and,
	\item $[f_C(X)+f_C(Y)]f_C(Z)=f_C(X)\cdot f_C(Z) +f_C(Y)\cdot f_C(Z)$
\end{enumerate}
\end{theorem}
\proof In the same way as Theorem \ref{thm.13}.
\qed


\subsection{Algebraic properties for the 'cross-ratio map sets' $\mathcal{R}^{D}_4$-set of points} $  $\\

We have that,
\[ 
\mathcal{R}^{D}_4=\left\{c_r(A,B;C,X) | \quad \forall X \in \ell^{OI} \right\}
\]

We mark elements of this set with $f_D(X)$ and, we have that
\[f_D(X)=c_r(A,B;C,X)=[(A-X)^{-1}(B-X)][(B-C)^{-1}(A-C)]
\]

We will see in the following Theorems, the properties that are fulfilled related to the addition and multiplication of points' in $\ell^{OI}-$line, so, we will see, the algebraic structures that are created, related to the addition and multiplication of this points in $\ell^{OI}-$line in Desargues affine plane.

We will not present proofs for some of the algebraic properties in this part, as they are the same as the proofs in the first part. We will present evidence only in cases where we have changes.

\begin{theorem} \label{thm.41}
Let's have three different-fixed points $A, B, C \in \ell^{OI}$  and this points are different for point $O$, in  $\ell^{OI}-$line. The transforms $f_D(\cdot)$ from $\mathcal{R}^{D}_4$, with addition of points, fulfills the properties: \\
(1) Associative,\\
(2) commutative and, \\
(3) it has a neutral element (zero element).
\end{theorem}
\proof  
\textbf{1.} Associative  and \textbf{2.} Commutative properties 
prove similar to Theorem \ref{thm.11}.

\textbf{3.} ('zero-element') $\exists Z \in \ell^{OI}, \forall X \in \ell^{OI}$, we have that:
\[ c_r(A,B;C,X)+c_r(A,B;C,Z)=c_r(A,B;C,X) \]
cross-ratio point are point of $\ell^{OI}$-line, and '$+$' is addition of points in $\ell^{OI}$-line, we have prove that the 'zero-element' for addition of points in $\ell^{OI}$-line, is point $O$, thus, have that,
\[ c_r(A,B;C,Z)=O \]
For point $O$ (zero element for addition of points), have
\[ c_r(A,B;C,Z)=O\]
so
\[ [(A-Z)^{-1}(B-Z)][(B-C)^{-1}(A-C)] =O\]

also we have that the skew field $(\ell^{OI},+, \cdot )$ do not have 'divisors of zero', therefore,
 \[ [(A-Z)^{-1}(B-Z)]=O \quad\text{or}\quad [(B-C)^{-1}(A-C)]=O \]
more,

 \[ (A-Z)^{-1}=O \quad\text{or}\quad (B-Z)=O , \]
or,
 \[ (B-C)^{-1}=O \quad\text{or}\quad (A-C)=O .\]

But, we suppose that the point $A \neq B$,  $A \neq C$ and $B \neq C$, therefore we have that,
\[
A-C \neq O \quad\text{and}\quad O\neq (B-C)^{-1} \neq \infty
\]
so have that,
\[ (A-Z)^{-1}=O \quad\text{or}\quad (B-Z)=O . \]
Hence, 
\[B-Z=O \Leftrightarrow Z=B
\]
so 'zero-element' is $f_D(B)$.
\qed

\begin{theorem}\label{thm.42}
Let's have three different points $A,C,D$ and different for point $O$, in  $\ell^{OI}-$line. The cross-ratio map set with multiplication of points in a line of Desargues affine plane $\left( \mathcal{R}^{B}_4, \cdot \right)$ is  Group.
\end{theorem}
\proof
\textbf{1.} Associative properties prove similar to Theorem \ref{thm.12}

\textbf{2.} ('unitary') $\exists E \in \ell^{OI}, \forall X \in \ell^{OI}$, we have that:
\[
f_D(X)\cdot f_D(E)=f_D(X)\]
so, 
\[
\{[(A-X)^{-1}(B-X)][(B-C)^{-1}(A-C)] \}
 \cdot \{[(A-E)^{-1}(B-E)][(B-C)^{-1}(A-C)] \}  \]
\[ =[(A-X)^{-1}(B-X)][(B-C)^{-1}(A-C)]\]
Both elements $f_D(X), f_D(E)$ are points in the line $\ell^{OI}$, we also know that the unitary element related to the multiplication of points in $\ell^{OI}-$ line is \emph{unique point} $I$, for this reason we have $f_D(E)=I$. We see that, for $E=C$, we have,
\[
\begin{aligned}
f_D(E)=f_D(C) \\
&=[(A-C)^{-1}(B-C)][(B-C)^{-1}(A-C)]\\
&=(A-C)^{-1}[(B-C)(B-C)^{-1}](A-C)\\
&=(A-C)^{-1}[I](A-C)\\
&=(A-C)^{-1}(A-C)\\
&=I.
\end{aligned}
\]
so unitary-element element is $f_D(C)$.

\textbf{3.} ('inverse') $\forall X \in \ell^{OI}, \exists Y \in \ell^{OI}$, we have that:
\[
f_D(X)\cdot f_D(Y)=I=f_D(C)\]
so,
\[
f_B(X)\cdot f_B(Y)=I\Rightarrow f_B(Y)=\{f_B(X) \}^{-1} \]
\[
\{[(A-X)^{-1}(B-X)][(B-C)^{-1}(A-C)] \}
 \cdot \{[(A-Y)^{-1}(B-Y)][(B-C)^{-1}(A-C)]\}=I, \]
hence,
\[
[(A-Y)^{-1}(B-Y)][(B-C)^{-1}(A-C)]=\{ [(A-X)^{-1}(B-X)][(B-C)^{-1}(A-C)] \}^{-1} \]
\[
=\{[(B-C)^{-1}(A-C)]\}^{-1}\{[(A-X)^{-1}(B-X)]\}^{-1} \]
\[f_D(Y) =[(A-C)^{-1}(B-C)][(B-X)^{-1}(A-X)] \]
Hence,
\[ \{f_D(X) \}^{-1}=c_r(A, B; X,C) \]
\qed

\begin{theorem} \label{thm.43}
In the cross-ratio-maps-set $\mathcal{R}^{D}_4=\left\{f_D(X)=c_r(A,B;C,X) | \quad \forall X \in \ell^{OI} \right\}$, for three elements $f_D(X)$, $f_D(Y)$, $f_D(Z)$ of its, have true, that equations
\begin{enumerate}
	\item $f_D(X)\cdot [f_D(Y)+f_D(Z)]=f_D(X)\cdot f_D(Y) +f_D(X)\cdot f_D(Z)$, and,
	\item $[f_D(X)+f_D(Y)]f_D(Z)=f_D(X)\cdot f_D(Z) +f_D(Y)\cdot f_D(Z)$
\end{enumerate}
\end{theorem}
\proof In the same way as Theorem \ref{thm.13}.
\qed

\section{Declarations}
\subsection*{Funding}
No Funding for this research.


\subsection*{Data availability statements:} This manuscript does not report data.
\subsection*{Conflict of Interest Statement}
There is no conflict of interest with any funder.


\begin{thebibliography}{99}

\bibitem{Artin1957GeometricAlgebra} Artin, E. {\it Geometric algebra}. Interscience Publishers, New York, NY,1957

\bibitem{Berger2010geometryRevealed} Berger, M. Geometry revealed. (Springer,2010), xvi+831 pp., ISBN: 978-3-540-70996-1, MR2724440

\bibitem{Berger2009geometry12} Berger, M. {\it Geometry. I, II}. Translated from the French by M. Cole and S. Levy. Berlin: Springer, 2009.

\bibitem{CoxterIG1969} Coxeter, H. {\it Introduction to geometry, 2nd Edition}. (John Wiley \& Sons, Inc.,1969), xvii+469 pp., MR0123930, MR0346644

\bibitem{Hartshorne1967Foundations} Hartshorne, R. {\it Foundations of projective geometry}. (New York: W.A. Benjamin, Inc. 1967. VII, 167 p. (1967).,1967)

\bibitem{Hilbert1959geometry} Hilbert, D. {\it The foundations of geometry}. (The Open Court Publishing Co.,1959), vii+143 pp., MR0116216

\bibitem{HugesPiper} Hughes, D. \& Piper, F. {\it Projective Planes}. Graduate Texts in Mathematics, Vol. 6. (Spnnger-Verlag,1973), x+291 pp., MR0333959

\bibitem{Kryftis2015thesis} Kryftis, A. {\it A constructive approach to affine and projective planes}. (University of Cambridge,2015), supervisor: M. Hyland, v+170pp.,arXiv 1601.04998v1 19 Jan. 2016

\bibitem{Pickert1973PlayfairAxiom} Pickert, G. {\it Affine Planes: An Example of Research on Geometric Structures}. {\em The Mathematical Gazette}. \textbf{57}, 278-291 (2004), MR0474017

\bibitem{Prazmowska2004DemoMathDesparguesAxiom} Prażmowska, M. A proof of the projective Desargues axiom in the Desarguesian affine plane. {\em Demonstratio Mathematica}. \textbf{37}, 921-924 (2004), MR2103894

\bibitem{Szmielew1981DesarguesAxiom} Szmielew, W. Od geometrii afinicznej do euklidesowej (Polish) [From affine geometry to Euclidean geometry] Rozwa?ania nad aksjomatyk? [An approach through axiomatics]. (Biblioteka Matematyczna [Mathematics Library],1981), 172 pp., ISBN: 83-01-01374-5, MR0664205



\bibitem{FilipiZakaJusufi} Filipi, K., Zaka, O. \& Jusufi, A. {\it The construction of a corp in the set of points in a line of Desargues affine plane}. {\em Matematicki Bilten}. \textbf{43}, 1-23 (2019), ISSN 0351-336X (print), ISSN 1857–9914 (online)

\bibitem{ZakaThesisPhd} Zaka, O. {\it Contribution to Reports of Some Algebraic Structures with Affine Plane Geometry and Applications}. (Polytechnic University of Tirana,Tirana, Albania,2016), supervisor: K. Filipi, vii+113pp.

\bibitem{ZakaCollineations} Zaka, O. {\it A description of collineations-groups of an affine plane}. {\em Libertas Mathematica (N.S.)}. \textbf{37}, 81-96 (2017), ISSN print: 0278 – 5307, ISSN online: 2182 – 567X, MR3828328

\bibitem{ZakaVertex} Zaka, O. {\it Three Vertex and Parallelograms in the Affine Plane: Similarity and Addition Abelian Groups of Similarly n-Vertexes in the Desargues Affine Plane}. {\em Mathematical Modelling And Applications}. \textbf{3}, 9-15 (2018), http://doi:10.11648/j.mma.20180301.12

\bibitem{ZakaDilauto} Zaka, O. {\it Dilations of line in itself as the automorphism of the skew-field constructed over in the same line in Desargues affine plane}. {\em Applied Mathematical Sciences}. \textbf{13}, 231-237 (2019). https://www.m-hikari.com/ams/ams-2019/ams-5-8-2019/p/zakaAMS5-8-2019.pdf

\bibitem{ZakaFilipi2016} Zaka, O. \& Filipi, K. {\it The transform of a line of Desargues affine plane in an additive group of its points}. {\em Int. J. Of Current Research}. \textbf{8}, 34983-34990 (2016). https://www.journalcra.com/sites/default/files/issue-pdf/16349.pdf

\bibitem{ZakaPeters2022DyckFreeGroup} Zaka, O., Peters, J.F. {\it DYCK FREE GROUP PRESENTATION OF POLYGON CYCLES IN THE RATIO OF COLLINEAR POINTS IN THE DESARGUES AFFINE PLANE}. J Math Sci 280, 605–630 (2024). https://doi.org/10.1007/s10958-024-06995-4

\bibitem{ZakaPeters2022InvariantPreserving} Orgest Zaka and James Peters. {\it Progress in Invariant and Preserving Transforms for the Ratio of Collinear Points in the Desargues Affine Plane Skew Field}.  Bol. Soc. Paran. Mat. (3s.) v. 2025 (43): 1-21.
https://doi.org/10.5269/bspm.68800

\bibitem{ZakaPeters2025CrosRatio} Orgest Zaka and James Peters. {\it Cross Ratio Geometry: Advances for Multiple Collinear Points in the Desargues Affine Plane}.  Bol. Soc. Paran. Mat. (3s.) v. 2025 (44): 1-16.


\bibitem{ZakaPeters2024}  Zaka, Orgest; Peters, James F. {\it Invariant and preserving transforms for cross ratio of 4-points in a line on Desargues affine plane}. Journal of Prime Research in Mathematics. 20, No. 2, 48-63 (2024). Zbl 07943368
https://jprm.sms.edu.pk/invariant-and-preserving-transforms-for-cross-ratio-of-4-points-in-a-line-on-desargues-affine-plane/

\bibitem{PetersZakaFundamentalGroup2023} Peters, J.F., Zaka, O. {\it Dyck fundamental group on arcwise-connected polygon cycles}. Afr. Mat. 34, 31 (2023). https://doi.org/10.1007/s13370-023-01067-3

\bibitem{ZakaMohammedEndo} Zaka, O. \& Mohammed, M. {\it The endomorphisms algebra of translations group and associative unitary ring of trace-preserving endomorphisms in affine plane}. {\em Proyecciones}. \textbf{39}, 821-834 (2020). https://doi.org/10.22199/issn.0717-6279-2020-04-0051

\bibitem{ZakaMohammedSF} Zaka, O. \& Mohammed, M. {\it Skew-field of trace-preserving endomorphisms, of translation group in affine plane.} {\em Proyecciones}. \textbf{39}, 835-850 (2020). https://doi.org/10.22199/issn.0717-6279-2020-04-0052

\bibitem{ZakaPetersIso} Zaka, O. \& Peters, J. {\it Isomorphic-dilations of the skew-fields constructed over parallel lines in the Desargues affine plane}. {\em Balkan J. Geom. Appl.}. \textbf{25}, 141-157 (2020), https://www.emis.de/journals/BJGA/v25n1/B25-1zk-ZBG89.pdf

\bibitem{ZakaPetersOrder} Zaka, O. \& Peters, J. {\it Ordered line and skew-fields in the Desargues affine plane}. {\em Balkan J. Geom. Appl.}. \textbf{26}, 141-156 (2021), https://www.emis.de/journals/BJGA/v26n1/B26-1zb-ZBP43.pdf

\bibitem{Cohn2008sf} Cohn, P. {\it Skew fields. Theory of general division rings}. Encycl. Math. Appl. ISSN: {0953-4806}, Volume {57}. Cambridge: Cambridge University Press (2008; Zbl 1144.16002)

\bibitem{Herstein1968NR} Herstein, I. {\it Topics in algebra, 2nd Ed}. (Xerox College Publishing,1975), pp. xi+388, MR0356988; first edition in 1964, MR0171801.

\bibitem{Rotman2015AMAlgebra} Rotman, J. {\it Advanced modern algebra. Part 1}. Providence, RI: American Mathematical Society (AMS), 2015.

\bibitem{Lam2001GTMalgebra} Lam, T. {\it A first course in noncommutative rings}. Graduate Texts in Mathematics, Vol.131. New York, NY: Springer,  (2001; Zbl 0980.16001).


\bibitem{Milne1911elementaryCross-Ratio} Milne, J. An elementary treatise on cross-ratio geometry, with historical notes.. (Cambridge: University Press. XXII u. 288 S. \(8^\circ\) (1911).,1911)


\end{thebibliography}

\end{document}